\documentclass[english]{llncs}

\usepackage{times}
\usepackage{amsfonts}
\usepackage{latexsym}
\usepackage{url}
\usepackage{epsfig}
\usepackage{amsfonts}
\usepackage{graphics}
\usepackage{graphicx}
\usepackage{llncsdoc} 

\newcommand{\st}{\,|\;}

\newcommand{\calU}{\mathcal{U}}
\newcommand{\calV}{\mathcal{V}}

\newcommand{\calG}{\mathcal{G}}
\newcommand{\calE}{\mathcal{E}}

\newcommand{\Ra}{\Rightarrow}
\newcommand{\La}{\Leftarrow}

\makeatother

\title{On graph equivalences preserved under extensions
}

\author{
Zbigniew Lonc\inst{1}
 \and
Miros{\l}aw Truszczy{\'n}ski\inst{2}}
\institute{%
Faculty of Mathematics and Information Science. Warsaw University of
Technology, 00-661 Warsaw, Poland,
\email{zblonc@mini.pw.edu.pl}
\and
Department of Computer Science, University of Kentucky,
Lexington, KY~40506, USA,
\email{mirek@cs.uky.edu}
}

\begin{document}
\maketitle

\begin{abstract}
Let $\cal G$ be the set of finite graphs whose vertices belong to
some fixed countable set, and let $\equiv$ be an equivalence
relation on $\cal G$. By the \emph{strengthening} of $\equiv$ we
mean an equivalence relation $\equiv_s$ such that $G\equiv_s H$,
where $G,H\in\cal G$, if for every $F\in \cal G$, $G\cup F\equiv
H\cup F$. The most important case that we study in this paper
concerns equivalence relations defined by graph properties. We write
$G\equiv^\Phi H$, where $\Phi$ is a graph property and $G,H\in\cal
G$,  if either both $G$ and $H$ have the property $\Phi$, or both do
not have it. We characterize the strengthening of the
relations $\equiv^\Phi$ for several graph properties $\Phi$. For
example, if $\Phi$ is the property of being a $k$-connected graph,
we find a polynomially verifiable (for $k$ fixed) condition that
characterizes the pairs of graphs equivalent with respect to
$\equiv_s^\Phi$. We obtain similar results when $\Phi$ is the
property of being $k$-colorable, edge $2$-colorable, hamiltonian, or
planar, and when $\Phi$ is the property of containing a subgraph
isomorphic to a fixed graph $H$. We also prove several general
theorems that provide conditions for $\equiv_s$ to be of some
specific form. For example, we find a necessary and sufficient
condition for the relation $\equiv_s$ to be the identity. Finally,
we make a few observations on the strengthening in a more general
case when $\cal G$ is the set of finite subsets of some countable
set.
\end{abstract}

\section{Introduction}

Equivalence relations partition their domains into classes of
equivalent objects --- objects indistinguishable with respect to
some characteristic. In the case of domains whose elements can
be combined, equivalence relations can be \emph{strengthened}. In
this paper, we introduce the concept of strengthening, motivate it,
and study it in the case of equivalence relations that arise in the
domain of graphs.

To illustrate what we have in mind, let us consider a set $X$ of
possible team members. Teams are finite subsets of $X$. We have some
equivalence relation on the set of teams, which groups in its
equivalence classes teams of the same value. Thus, given two
equivalent teams, say $A,B\subseteq X$, we could use any of them
without compromising the quality. But there is more to it. Let us
consider a team $C$, such that $A$ is its sub-team, that is,
$A\subseteq C$. Let us also suppose that for one reason or another
we are unable to keep all members of $A$ in $C$. If we need the
``functionality'' of $A$ in $C$, we might want to replace $A$ with
its equivalent $B$ by forming the team $C'=B\cup(C\setminus A)$.
After all, $A$ and $B$ are equivalent. But this is a reasonable
solution only if by doing so, we do not change the quality of the
overall team, that is, if $C$ and $C'$ are equivalent, too. And, in
general, it is not guaranteed.

Let us observe that in our example $C=A\cup(C\setminus(A\cup B))$
and $C' =B\cup(C\setminus (A\cup B)$, that is, they are extensions
of $A$ and $B$, respectively, with the same set, $(C\setminus(A\cup
B))$. This suggests that we might call teams $A$ and $B$
\emph{strongly equivalent} (with respect to the original equivalence
relation) if for every finite set $D$, $A\cup D$ and $B\cup D$ are
equivalent. Clearly, if $A$ and $B$ are strongly equivalent, then
any two teams obtained by extending $A$ and $B$ with the same
additional members are equivalent! Thus, the relation of strong
equivalence, the ``strengthening'' of the original one, is precisely
what we need when we consider teams not as individual entities but
as potential sub-teams in bigger groups.

To the best of our knowledge, the concept of strong equivalence has
emerged so far only in the area of logic
programming \cite{lpv01,tu03,efw04,Woltran07}.
Researchers argued there that it underlies the notion of
a module of a program, and is essential to modular program development.
In this paper we study the strengthening of an equivalence relation in the
domain of graphs. As a result, we obtain a new class of graph-theoretic
problems. Importantly, when applied to specific properties, for instance,
to the graph connectivity, the notion of strengthening has interesting
practical implication and does give rise to non-trivial arguments and
characterizations.

Let us consider the following scenario. In the context of networks,
which we typically represent as graphs, the concept of their
connectivity is of paramount importance (cf. Colbourn \cite{C}). Let
us define graphs $G$ and $H$ to be equivalent if both are
$k$-connected or if neither of them is. With time networks grow and
get embedded into bigger networks. The key question is: are the
graphs $G$ and $H$ interchangeable, in the sense that the networks
obtained by identical extensions of $G$ and $H$ are equivalent with
respect to $k$-connectivity?

\begin{figure}
\centerline{\includegraphics[scale=0.40]{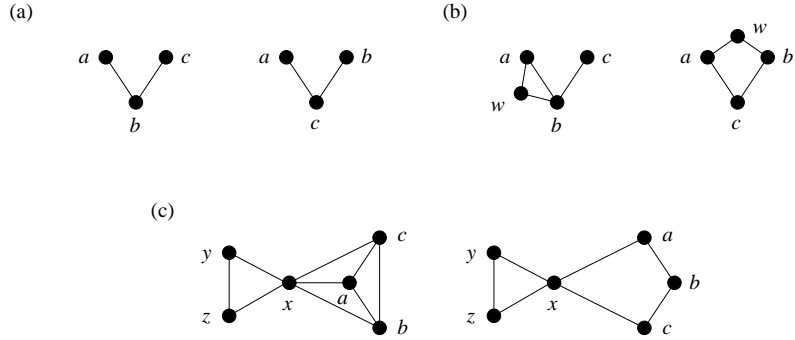}}
\caption{(a) Graphs that are not strongly equivalent with respect to
2-connectivity; (b) Extensions of the graphs from (a) showing that graphs
in (a) are not strongly equivalent with respect to 2-connectivity; and (c)
Two graphs that are strongly equivalent with respect to
2-connectivity.}
\label{mot}
\end{figure}

For example, neither of the two graphs in Figure \ref{mot}(a) is
2-connected and so, they are equivalent with respect to
2-connectivity. They are not, however, strongly equivalent with
respect 2-connectivity. Indeed, the graphs obtained by extending
them with two edges $aw$ and $bw$, shown in Figure \ref{mot}(b) are
not 2-connectivity equivalent --- one of them is 2-connected and the
other one is not! On the other hand, one can verify directly from the
definition that the graphs shown in Figure \ref{mot}(c) are strongly
equivalent with respect to 2-connectivity. Later in the paper, we provide
a characterization that allows us to decide the question of strong
equivalence with respect to connectivity.

Our paper is organized as follows. While most of our results concern
strengthening of equivalence relations on graphs, we start by introducing
the concept of the strengthening of an equivalence relation in a more
general setting of the domain of finite subsets of a set. We derive there
several basic properties of the notion, which we use later in the paper.
In particular, for a class of equivalence relations defined in terms of
properties of objects --- such equivalence relations are of primary
interest in our study --- we characterize those relations that are equal
to their strengthenings.

The following sections are concerned with equivalence relations on graphs
defined in terms of graph properties, a primary subject of interest to us.
Narrowing down the focus of our study to graphs allows us to obtain stronger
and more interesting results. In particular, we characterize relations whose
strengthening is the identity relation, and those whose strengthening
defines one large equivalence class, with all other classes being
singletons. We apply these general characterizations to obtain
descriptions of strong equivalence with respect to several concrete
graph-theoretic properties including possessing hamiltonian cycles and
being planar. Main results of the paper, are contained in the two sections
that follow. They concern graph-theoretic properties, which do not fall
under the scope of our general results. Specifically, we deal there with
vertex and edge colorings, and with $k$-connectivity. The characterizations
we obtain are non-trivial and show that the idea of strengthening gives rise
to challenging problems that often (as in the case of strengthening
equivalence with respect to $k$-connectivity) have interesting motivation
and are of potential practical interest.

\section{The Problem and General Observations}

We fix an infinite countable set $\calE$ and denote by $\calG$ the
set of finite subsets of $\calE$.

\begin{definition}
Let $\equiv$ be an equivalence relation on $\calG$. We say that sets
$G,H\in\calG$ are \emph{strongly equivalent} with respect to $\equiv$,
denoted by $G\equiv_s H$, if for every set $F\in\calG$, $G\cup F \equiv
H\cup F$. We call $\equiv_s$ the \emph{strengthening} of $\equiv$.
\end{definition}

While most of our results concern the case when $\calE$ is a set of
edges over some infinite countable set of vertices $\calV$, in this
section we impose no structure on $\calE$ and prove several basic
general properties of the concept of strong equivalence.

\begin{proposition}
\label{prop:1}
Let $\equiv$ be an equivalence relation on sets in $\calG$. Then:
\begin{enumerate}
\item the relation $\equiv_s$ is an equivalence relation
\item for every sets $G,H\in\calG$, $G\equiv_s H$ implies $G\equiv H$
(that is, $\equiv_s\ \subseteq\ \equiv$)
\item for every sets $G,H,F\in\calG$, $G\equiv_s H$ implies $G\cup F\equiv_s
H\cup F$.
\end{enumerate}
\end{proposition}
Proof: (1) All three properties of reflexivity, symmetry and transitivity
are easy to check. For instance, let us assume that for some three sets
$D,G,H\in\calG$, $D\equiv_s G$ and $G\equiv_s H$. Let $F\in\calG$. By the
definition, $D\cup F\equiv G\cup F$ and $G\cup F\equiv H\cup F$. By the
transitivity of $\equiv$, $D\cup F\equiv H\cup F$. Since $F$ was an
arbitrary element of $\calG$, $D\equiv_s H$.

\smallskip
\noindent
(2) By the definition of $\equiv_s$, for every set $F\in\calG$, $G\cup
F\equiv H\cup F$. In particular, if $F=\emptyset$, we get that $G\equiv H$.

\smallskip
\noindent
(3) For every set $F'\in\calG$, $F\cup F'\in\calG$. Since, $G\equiv_s H$,
$(G\cup F)\cup F' = G\cup(F\cup F') \equiv_s H\cup(F\cup F') = (H\cup F)\cup
F'$. Thus, the claim follows.
\hfill$\Box$

\begin{proposition}
\label{prop:2}
Let $\approx$ and $\equiv$ be equivalence relations on $\calG$. Then:
\begin{enumerate}
\item if $\approx\ \subseteq\ \equiv$, then $\approx_s\ \subseteq\ \equiv_s$
\item if $\approx_s\ =\ \equiv_s$, then $(\approx\cap\equiv)_s\ =\ \approx_s$.
\end{enumerate}
\end{proposition}
Proof: Arguments for each of the assertions are simple. As an example,
we prove (2) here. By (1), it suffices to show that $\approx_s\subseteq
(\approx \cap\equiv)_s$. Thus, let us consider sets $G,H\in\calG$ such that
$G\approx_s H$. Let $F\in\calG$. Clearly, $G\cup F\approx H\cup F$. Moreover,
by the assumption, $G\equiv_s H$. Thus, $G\cup F\equiv H\cup F$, as well.
It follows that $G\cup F (\approx\cap\equiv) H\cup F$. As $F$ is arbitrary,
$G(\approx\cap\equiv)_s H$ follows. \hfill$\Box$

\begin{corollary}
\label{cor:fixed}
Let $\equiv$ be an equivalence relation on sets in $\calG$. Then,
$(\equiv_s)_s = \equiv_s$. Moreover, for every equivalence relation
$\approx$ on $\calG$ such that $\approx_s = \equiv_s$, $\equiv_s \subseteq
\approx$.
\end{corollary}
Proof: By Proposition \ref{prop:1}(2), $\equiv_s\subseteq\equiv$. Thus,
by Proposition \ref{prop:2}(1), $(\equiv_s)_s\subseteq\equiv_s$. Conversely,
let us consider sets $G,H\in\calG$ such that $G\equiv_s H$ and let
$F,F'\in\calG$. Since $G\equiv_s H$, $G\cup(F\cup F')\equiv H\cup(F\cup F')$.
Thus, $(G\cup F)\cup F'\equiv (H\cup F)\cup F'$. As $F'$ is arbitrary, $G\cup
F\equiv_s H\cup F$ follows. Consequently, as $F$ is arbitrary, too,
$G(\equiv_s)_s H$.

To prove the ``moreover'' part of the assertion, we note that $\approx_s
\subseteq\approx$. Thus, $\equiv_s\subseteq\approx$, as needed. \hfill$\Box$

\smallskip
Corollary \ref{cor:fixed} states, in particular, that for every equivalence
relation $\equiv$ on $\calG$, $\equiv_s$ is the \emph{most precise} among all
equivalence relations $\approx$ such that $\approx_s=\equiv_s$.

Most of our results concern equivalence relations defined in terms of
functions assigning to sets in $\calG$ collections of certain objects.
Let $\calU$ be a set and let $f:\calG\rightarrow2^\calU$. For sets $G$ and $H$,
we define:
\begin{enumerate}
\item $G\cong^f H$ if $f(G)=f(H)$, and
\item $G\equiv^f H$ if $f(G)=f(H)=\emptyset$, or $f(G)\not=\emptyset$ and
$f(H)\not=\emptyset$.
\end{enumerate}
Obviously, $\cong^f\ \subseteq\ \equiv^f$. Thus, by our earlier results,
$\cong^f_s\ \subseteq\ \equiv^f_s$, $\equiv^f_s\ \subseteq\ \equiv^f$, and
$\cong^f_s\ \subseteq\ \cong^f$.

\emph{Properties} of elements of $\calG$ (formally, subsets of
$\calG$) give rise to a special class of equivalence relations of
the latter type. Namely, given a property $\Phi\subseteq\calG$, we
define $\calU=\{0\}$ and set $f_\Phi(G)=\{0\}$ if and only if
$G\in\Phi$ (otherwise, $f_\Phi(G) =\emptyset$). Clearly,
$G\equiv^{f_\Phi}H$ if and only if both $G$ and $H$ have $\Phi$
($G,H\in\Phi$), or if neither $G$ nor $H$ does ($G\notin\Phi$ and
$H\notin\Phi$). To simplify the notation, we always write
$\equiv^\Phi$ instead of $\equiv^{f_\Phi}$. By $\overline{\Phi}$ we
denote the property $\calG - \Phi$.

In the remainder of this section we present a general result
concerning the relation $\equiv^\Phi$ that does not require any
additional structure of subsets of $\calG$. It characterizes those
properties $\Phi\subseteq \calG$, for which
$\equiv^\Phi=\equiv^\Phi_s$ (the strengthening does not change the
equivalence relation). The remainder of the paper is concerned with
the relations $\cong^f$ and $\equiv^f$ (including relations
$\equiv^\Phi$) and their strengthenings in the case when $\calG$
consists of graphs. In several places, we will consider properties
that are monotone. Formally, a property $\Phi\subseteq\calG$ is
\emph{monotone} if for every $G,H\in\calG$, $G\subseteq H$ and
$G\in\Phi$ imply $H\in\Phi$.

\begin{lemma}
\label{lem:equal} Let $\Phi\subseteq\calG$ be a property such that
$\emptyset\notin\Phi$. Then, $\equiv^\Phi = \equiv^\Phi_s$ if and
only if there is $X\subseteq \calE$ such that $\Phi=\{G\in\calG\st
G\cap X\not=\emptyset\}$.
\end{lemma}
Proof: $(\La)$ Let us assume that there is $X\subseteq\calE$ such
that $\Phi= \{G\in\calG\st G\cap X\not=\emptyset\}$. To prove that
$\equiv^\Phi = \equiv^\Phi_s$, it suffices to show that
$G\equiv^\Phi H$ implies $G \equiv^\Phi_s H$ (the converse
implication follows by Proposition \ref{prop:1}(2)).

Thus, let $G\equiv^\Phi H$. First, let us assume that $G,H\in\Phi$.
It follows that $G\cap X\not=\emptyset$ and $H\cap X\not=\emptyset$.
Consequently, for every set $F$, $(G\cup F)\cap X\not=\emptyset$ and
$(H\cup F)\cap X\not=\emptyset$. Thus, $G\cup F,H\cup F\in\Phi$ and
so, $G\cup F\equiv^\Phi H\cup F$. It follows that $G\equiv^\Phi_s
H$.

The only remaining possibility is that $G,H\notin\Phi$. Since $G\cap X=H
\cap X = \emptyset$, for every graph $F$, $(G\cup F)\cap X\not=\emptyset$
if and only if $(H\cup F)\cap X\not=\emptyset$. That is, for every graph
$F$, $G\cup F\equiv^\Phi H\cup F$ and so, $G\equiv^\Phi_s H$ in this case,
as well.

\smallskip
\noindent
$(\Ra)$ First, we prove that $\Phi$ is monotone.
Let $G\subseteq H$ and
let $G\in\Phi$. Let us assume that $H\notin\Phi$. It follows that
$\emptyset\equiv^\Phi H$ and so, by the assumption, $\emptyset\equiv^\Phi_s
H$. Thus, $G=\emptyset\cup G\equiv^\Phi H\cup G=H$. Since $H\notin\Phi$,
$G\notin\Phi$, a contradiction. Thus, $H\in\Phi$.

We define $X=\{e\in\calE\st \{e\}\in\Phi\}$. We will show that
$G\in\Phi$ if and only if $G\cap X\not=\emptyset$. If $G\cap X\not=
\emptyset$, then there is $e\in X$ such that $\{e\}\subseteq G$. Since
$\{e\}\in\Phi$, by the monotonicity of $\Phi$ it follows that $G\in\Phi$.

Conversely, let us assume that $G\in\Phi$. Let $G'\subseteq G$ be a
maximal subset of $G$ such that $G'\notin\Phi$. Such $G'$ exists as
$\emptyset \notin\Phi$. Moreover, since $G\in\Phi$, $G'\not=G$. It
follows that there is $e\in G\setminus G'$. We have
$\emptyset\equiv^\Phi G'$ as neither set has property $\Phi$. By the
assumption, $\emptyset\equiv^\Phi_s G'$. Thus, $\{e\}=\emptyset\cup
\{e\} \equiv^\Phi G'\cup\{e\}$. By the maximality of $G'$,
$G'\cup\{e\}\in\Phi$. Thus, $\{e\}\in \Phi$. Consequently, $e\in X$
and $G\cap X\not=\emptyset$. \hfill$\Box$

\begin{theorem}
Let $\Phi$ be a property. Then, $\equiv^\Phi = \equiv^\Phi_s$ if and only
if there is $X\subseteq \calE$ such that $\Phi=\{G\in\calG\st G\cap X\not=
\emptyset\}$ or $\Phi=\{G\in\calG\st G\subseteq X\}$.
\label{thm:equal}
\end{theorem}
Proof: $(\Ra)$ Let us assume that $\emptyset\not\in\Phi$. Then, Lemma
\ref{lem:equal} implies that there is a set $X$ such that
$\Phi=\{G\in\calG\st G\cap X\not=\emptyset\}$. If $\emptyset\in\Phi$,
then $\emptyset\notin\overline{\Phi}$. Since $\equiv^\Phi=
\equiv^{\overline{\Phi}}$, $\equiv^\Phi_s= \equiv^{\overline{\Phi}}_s$. Thus,
$\equiv^{\overline{\Phi}} = \equiv^{\overline{\Phi}}_s$. By Lemma
\ref{lem:equal},
there is a set $Y$ such that $\overline{\Phi}=\{G\in\calG\st G\cap Y\not=
\emptyset\}$. Setting $X=\calE\setminus Y$, we obtain that $\Phi=\{G\in\calG\st G\subseteq X\}$.

\smallskip
\noindent
$(\La)$ If $\Phi=\{G\in\calG\st G\cap X\not= \emptyset\}$, then $\emptyset
\notin\Phi$ and the result follows from Lemma \ref{lem:equal}. Thus, let
us assume that $\Phi=\{G\in\calG\st G\subseteq X\}$. We have then that
$\overline{\Phi}=\{G\in\calG\st G\cap Y\not=\emptyset\}$, where $Y=\calE
\setminus X$. Moreover, $\emptyset\notin\overline{\Phi}$. Thus, by Lemma
\ref{lem:equal}, $\equiv^{\overline{\Phi}}_s=\equiv^{\overline{\Phi}}$.
By our earlier observations, $\equiv^{\Phi}_s=\equiv^{\Phi}$.
\hfill$\Box$

\section{Strengthening of equivalence relations on graphs}

 From now on we focus on the special case of equivalence relations on
graphs. That is, we assume a fixed infinite countable set $\calV$ of
\emph{vertices} and define $\cal E$ to be the set of all unordered
pairs of two distinct elements from $\cal V$ (the set of all
\emph{edges} on $\cal V$). Thus, the elements of $\calG$ (finite
subsets of $\calE$) can now be regarded as graphs with the vertex
set implicitly determined by the set of edges. For a graph $G$, we
denote by $V(G)$ the set of vertices of $G$, that is, the set of all
endvertices of edges in $G$. From this convention, it follows that
we consider only graphs with no isolated vertices. We understand the
union of graphs as the set-theoretic union of sets. Given a graph
$G$ and edges $e\notin G$ and $f\in G$, we often write $G+e$ and
$G-f$ for $G\cup\{e\}$ and $G-\{f\}$, respectively. We refer the
reader to Diestel \cite{D} for definitions of all graph theoretic
concepts not defined in this paper.

Our first result fully characterizes equivalence relations of $\calG$
whose strengthening is the identity relation.

\begin{theorem}
\label{thm:id}
Let $\equiv$ be an equivalence relation on $\calG$. Then,
$\equiv_s$ is the identity relation on $\calG$ if and only
if for every complete graph $K\in\calG$ and for every $e\in K$,
$K\not\equiv_s K-e$.
\end{theorem}
Proof: $(\Ra)$ Let $K$ be a complete graph in $\calG$
and let $e\in K$. Since $K\not=K-e$, $K\not\equiv_s K-e$, as
required.

\smallskip
\noindent
$(\La)$ Let $G,H\in\calG$. Clearly, if $G=H$ then $G\equiv_s H$ (as
$\equiv_s$ is reflexive). Conversely, let $G\equiv_s H$. Let us assume
that $G\not=H$. Without loss of generality, we may assume that there is
an edge $e\in G\setminus H$. Let $K$ be the complete graph on the set
of vertices $V(G\cup H)$. Since $G\equiv_s H$, $G\cup(K-e) \equiv_s H\cup
(K-e)$ (cf. Proposition \ref{prop:1}(3)). Moreover, $G\cup(K-e)=K$ (as
$G\subseteq K$ and $e\in G$) and $H\cup(K-e)=K-e$ (as $H\subseteq K$ and
$e\notin H$). Consequently, $K\equiv_s K-e$, a contradiction. It follows
that $G=H$. \hfill$\Box$

\begin{remark} If we assume (as we did in the previous section) that
$\calG$ is simply a collection of all finite subsets of an arbitrary
infinite
countable set $\calE$, we could prove the following result (by essentially
the same method we used in Theorem \ref{thm:id}): Let $\equiv$ be an
equivalence relation on $\calG$. Then, $\equiv_s$ is the identity relation
on $\calG$ if and only if for every finite $S\in\calG$ and every $e\in S$,
$S\not\equiv_s S-e$. Applying this result to the case when $\calE$ is the
set of edges over some infinite countable set of vertices, and $\calG$ 
is a set of graphs built of edges in $\calE$ gives a weaker characterization
than the one we obtained, as its condition becomes ``for every graph $S$ 
and every $e\in S$, $S\not \equiv_s S-e$,'' while the condition in Theorem
\ref{thm:id} is
restricted to complete graphs only.
\end{remark}

We will now illustrate the applicability of this result. Let
$G,H\in\calG$ and let us define $G\equiv^{hc}H$ if and only if $G$
and $H$ either both have or both do not have a hamiltonian cycle.
Theorem \ref{thm:id} implies that the relation $\equiv^{hc}_s$ is
the identity relation. In other words, for every two distinct graphs
$G$ and $H$, there is a graph $F$ such that exactly one of the
graphs $G\cup F$ and $H\cup F$ is hamiltonian.

\begin{theorem}
\label{thm:hc}
Let $G,H\in\calG$. Then $G\equiv^{hc}_s H$ if and only if $G=H$.
\end{theorem}
Proof: $(\Ra)$ Let $K$ be a finite complete graph from $\calG$ and
let $e\in K$. If $|K|=1$, then $K=\{e\}$ and $K-e=\emptyset$. Let
$e=uv$ and let $w$ be a vertex in $\calV$ different from $u$ and
$v$. We define $F= \{uw,vw\}$. Clearly, $K\cup F$ has a hamiltonian
cycle and $(K-e)\cup F=F$ does not have one. Thus,
$K\not\equiv^{hc}_s K-e$. Next, let us assume that $|K|=3$. Then,
$K$ has a hamiltonian cycle and $K-e$ does not. Thus,
$K\not\equiv^{hc} K-e$ and, consequently, $K\not\equiv^{hc}_s K-e$.

\begin{figure}
\centerline{\includegraphics[scale=0.40]{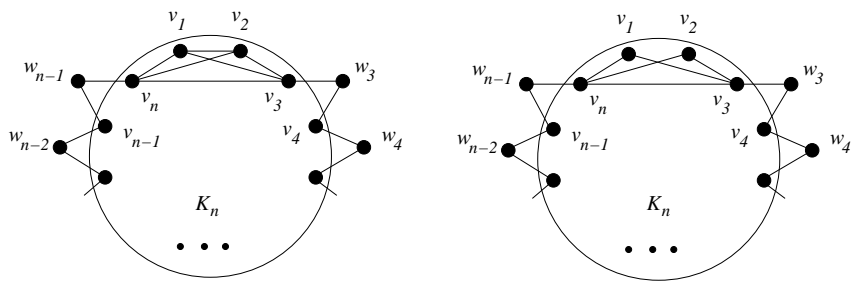}}
\caption{Graphs $K\cup F$ and $(K-e)\cup F$}
\label{fig1}
\end{figure}

Finally, let us assume that $|K|\geq 6$ (that is, $K$ is a finite
complete graph on at least 4 vertices). Let $v_1,\ldots, v_n$, where
$n\geq 4$, be the vertices of $K$ and let $e=v_1v_2$. We select
fresh vertices from $\calV$, say $w_3,\ldots, w_{n-1}$. We set
$F=\{v_iw_i\st i=3,\ldots,n-1\}\cup\{w_iv_{i+1}\st
i=3,\ldots,n-1\}$. It is clear that $K\cup F$ has a hamiltonian
cycle. However, $(K-e)\cup F$ does not have one! Indeed, any such
cycle would have to contain all edges in $F$, and that set cannot be
extended to a hamiltonian cycle in $(K-e)\cup F$ (cf. Figure
\ref{fig1}). Thus, also in this case $K\not\equiv^{hc}_s K-e$ and,
by Theorem \ref{thm:id}, $\equiv^{hc}_s$ is the identity relation.

\smallskip
\noindent
$(\La)$ This implication is evident. If $G=H$ then, clearly, $G\equiv^{hc}_s
H$. \hfill$\Box$

We note that the strengthening of a related (and more precise) equivalence
relation $\cong^{hc}$, where $hc$ is the function that assigns to a graph
the set of its hamiltonian cycles, is also the identity relation.

Next, we turn our attention to other equivalence relations on graphs
determined by properties of graphs (subsets of $\calG$).

We say that a property $\Phi\subseteq\calG$ is \emph{strong} if for
every $G,H\in\calG$, $G\equiv^{\Phi}_s H$ if and only if
$G,H\in\Phi$ or $G=H$. In other words, a property $\Phi$ is strong,
if the strong equivalence with respect to $\equiv^{\Phi}$ does not
``break up'' the equivalence class $\Phi$ of the relation $\equiv^{\Phi}$
(does not distinguish between
any graphs with the property $\Phi$) but, in the same time, breaks
the other equivalence class into singletons. 

We will now characterize properties $\Phi$ that are strong.

\begin{theorem}
\label{thm:exchange}
Let $\Phi$ be a property of graphs in $\calG$ (a subset of $\calG$).
Then, $\Phi$ is strong if and only if $\Phi$ is monotone and for every
graph $G\notin\Phi$, and every edge $e\in G$, $G\not\equiv^\Phi_s G-e$.
\end{theorem}
Proof: $(\La)$ We need to show that for every $G,H\in\calG$,
$G\equiv^{\Phi}_s H$ if and only if $G,H\in\Phi$ or $G=H$. If
$G,H\in\Phi$ then, by the monotonicity of $\Phi$, for every $F\in
\calG$, $G\cup F\in\Phi$ and $H\cup F\in\Phi$. Thus, $G\cup F\equiv^{\Phi}
H\cup F$. Since $F$ is arbitrary, $G\equiv^{\Phi}_s H$. If $G=H$ then
$G\equiv^{\Phi}_s H$ is evident.

Conversely, let $G\equiv^{\Phi}_s H$. Then $G\equiv^{\Phi} H$ and so,
either both $G$ and $H$ have $\Phi$, or neither $G$ nor $H$ has $\Phi$.
In the first case, there is nothing left to prove. Thus, let us assume
that $G,H\notin\Phi$.

Let $e\in G\setminus H$. We define $G'=G\cup H$. Since $G\equiv^{\Phi}_s
H$, $G'=G\cup H\equiv^{\Phi} H\cup H=H$. Since $H\notin\Phi$, $G'\notin
\Phi$. Let $G''$ be a maximal graph such that $V(G'')=V(G')$, $G'\subseteq
G''$ and $G''\notin\Phi$ (since $G'\notin\Phi$, such a graph exists). We
observe that $G''=G\cup (G''-e)$ and $G''-e=H\cup(G''-e)$. Thus, $G''
\equiv^{\Phi}_s (G''-e)$ (cf. Proposition \ref{prop:1}(3)), a contradiction.
It follows that $G\subseteq H$. By symmetry, $H\subseteq G$ and so, $G=H$.

\smallskip
\noindent
$(\Ra)$ Let us assume that there is a graph $G\notin\Phi$ such that for
every supergraph $G'$ of $G$, $G'\notin\Phi$. Let $H$ be any \emph{proper}
supergraph of $G$ (clearly, $G$ has proper supergraphs). and let $F$ be
any graph. Then $G\cup F$ and $H\cup F$ are both supergraphs of $G$. It
follows that $G\cup F\notin\Phi$ and $H \cup F\notin\Phi$. Consequently,
$G\cup F\equiv^\Phi H\cup F$. Since $F$ is an arbitrary graph,
$G\equiv^\Phi_s H$. However, $\Phi$ is strong and so the equivalence
class of $G$ under $\equiv^\Phi_s$ consists of $G$ only (as $G\notin\Phi$).
Thus, $G=H$, a contradiction.
It follows that for every graph $G\notin\Phi$, there is a supergraph
$G'$ of $G$ such that $G'\in\Phi$.

Let us now assume that $\Phi$ is not
monotone. Then, there are graphs $G$ and $H$ such that $G\subseteq H$,
$G\in\Phi$, and $H\notin\Phi$. Let $H'$ be a supergraph of $H$ such
that $H'\in\Phi$. It follows that $G\equiv^\Phi H'$ and, as $\Phi$ is
strong, $G\equiv^\Phi_s H'$. Thus, $H=G\cup H \equiv^\Phi H'\cup H= H'$,
a contradiction (as $H\notin\Phi$ and $H'\in\Phi$). It follows that
$\Phi$ is monotone.

Finally, let $G\notin\Phi$ and $e\in G$. Since $\Phi$ is strong, the
equivalence class of $G$ under $\equiv^\Phi_s$ consists of $G$ only.
Consequently, $G\not\equiv^\Phi_s G-e$.
\hfill$\Box$

To illustrate the scope of applicability of this result, we will
consider now several graph-theoretic properties. We start with the
property of non-planarity, that is, the set of all graphs that are
not planar.

\begin{theorem}
\label{thm:non-pl}
The property of non-planarity is strong.
\end{theorem}
Proof: Let $\Phi$ denote the property of non-planarity. It is clear that
$\Phi$ is monotone. Thanks to Theorem \ref{thm:exchange}, to complete
the proof it suffices to show that for every graph $G\notin\Phi$ and every
edge $e\in G$, $G\not\equiv^\Phi_s G-e$.

Thus, let $G\notin\Phi$, that is, let $G$ be a planar graph. Let $e\in
G$. We will denote by $x$ and $y$ the endvertices of $e$. First, let us
assume that $G=\{e\}$. Let $K$ be a complete graph on 5 vertices that
contains $e$. Clearly, $\emptyset \cup (K-e)$ is planar and $\{e\}\cup
(K-e)=K$ is not. Thus, $G\not\equiv^\Phi_s G-e$.

 From now on, we will assume that $G$ has at least three vertices.
Let $G'$ be a maximal planar supergraph of $G$ such that
$V(G)=V(G')$. We will now fix a particular planar embedding of $G'$,
and assume that $e$ belongs to the outerface (such an embedding
exists). With some abuse of terminology, we will refer also to this
embedding as $G'$. We observe that by the maximality of $G'$, every
face in $G'$ is a triangle. Since $G\subseteq G'$, to prove that
$G\not\equiv^\Phi_s G-e$, it suffices to show that
$G'\not\equiv^\Phi_s G'-e$ (by Proposition \ref{prop:1}(3)).

\smallskip
\noindent 1. $|V(G')|=3$. Then $G'$ is a triangle. Let $x,y$, and
$z$ be the vertices of $G'$ and let (as before $e=xy$). Let $v$ and
$w$ be two new vertices and $F=\{vx,vy,vz,wx,wy,wz,vw\}$. Then
$G'\cup F=K$, where $K$ is a complete graph on 5 vertices. Clearly,
$K$ is not planar. On the other hand, $(G'-e) \cup F=K-e$ is planar.
Thus, $G'\cup F \not\equiv^\Phi(G'-e)\cup F$ and so,
$G'\not\equiv^\Phi_s G'-e$.

\begin{figure}
\centerline{\includegraphics[scale=0.40]{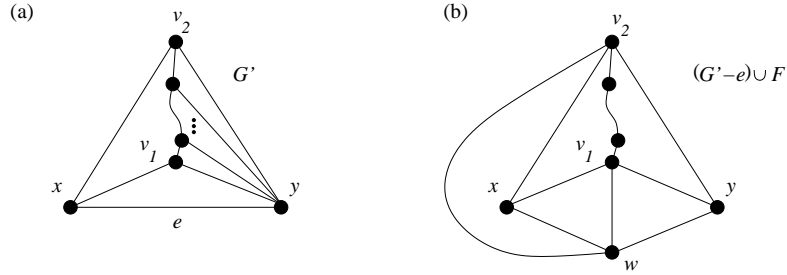}}
\caption{(a) Graph $G'$; (b) Graph $G'\cup F$.}
\label{pl}
\end{figure}

\smallskip
\noindent 2. $|V(G')|\geq 4$. There are two distinct faces, say
$F_1$ and $F_2$ in $G'$, sharing the edge $e$. Both faces are
triangles and, without loss of generality, we assume that $F_2$ is
the outerface. Let us assume, as before, that $e=xy$ and let $v_1$
(respectively, $v_2$) be the third vertex of the face $F_1$ (respectively,
$F_2$). We note that there is a path from $v_1$ to $v_2$ in $G'$
that does not contain $x$ nor $y$. Indeed, every two edges incident
to $y$ and consecutive in the embedding of $G'$ are connected with
an edge, as all faces in $G'$ are triangles (cf. Figure
\ref{pl}(a)).

Let $w$ be a new vertex and $F=\{wx,wy,wv_1,wv_2\}$. Then $(G'-e)
\cup F$ is planar (cf. Figure \ref{pl}(b)). On the other hand, the
graph $G'\cup F$ contains a subgraph homomorphic to the complete
graph on five vertices and so, $G\cup F$ is not planar. Thus,
$G'\cup F\not\equiv^\Phi(G'-e) \cup F$ and, consequently,
$G'\not\equiv^\Phi_s G'-e$. \hfill$\Box$

\begin{corollary}
Let $\equiv^{pl}$ be the equivalence relation such that for every two
graphs $G$ and $H$, $G\equiv^{pl}H$ if both $G$ and $H$ are planar or
both $G$ and $H$ are non-planar. Then, for every two graphs $G$ and $H$,
$G\equiv^{pl}_s H$ if and only if both $G$ and $H$ are non-planar or
$G=H$.
\end{corollary}
Proof: The relation $\equiv^{pl} = \equiv^\Phi$, where $\Phi$ is the
non-planarity property. Since the relation $\Phi$ is strong (by Theorem
\ref{thm:non-pl}), the assertion follows. \hfil$\Box$

Theorem \ref{thm:exchange} applies to many graph-theoretic properties concerned
with the containment of particular subgraphs. We will present several
such properties below.

\begin{theorem}
\label{thm:pm}
The property $\Phi_H$ consisting of all graphs containing a subgraph
isomorphic to $H$ is strong in each of the following cases:
\begin{enumerate}
\item $H$ is a star
\item $H$ is a cycle
\item $H$ is a 2-connected graph such that for every 2-element cutset $\{x,y\}$,
$xy$ is an edge of $H$
\item $H$ is a 3-connected graph
\item $H$ is a complete graph
\end{enumerate}
\end{theorem}
Proof: In each case the property is monotone. Thus, we only need to show
that for every $G\notin\Phi_H$ and every edge $e\in G$, there is a graph
$F$ such that $G\cup F\not\equiv^{\Phi_H} (G-e)\cup F$. Below we assume
that $e=xy$.

\smallskip
\noindent (1) If $H=\{e\}$, for some $e\in\calE$, and
$G\notin\Phi_H$, then $G=\emptyset$ and so the required property
holds vacuously. Thus, let us assume that $H$ consists of $k\geq 2$
edges. Since $G\notin\Phi_H$, $deg_G(x)<k$. Let $F$ be a star with
$k-deg_G(x)$ edges all incident to $x$ and with the other end not
in $G$. Clearly, $G\cup F$ contains a star with $k$ edges. On the
other hand, $(G-e)\cup F$ does not.

\smallskip
\noindent
(2) Let $H$ be a cycle with $k$ edges. We define $F$ to be a path with
$k-1$ edges, with endvertices $x$ and $y$, and with all intermediate vertices
not in $G$. Clearly, $G\cup F$ has a cycle of length $k$ and $(G-e)\cup F$
does not.

\smallskip
\noindent (3) Let $F'$ be a graph isomorphic to $H$, such that
$V(F')\cap V(G)=\{x,y\}$ and $xy$ is an edge of $F'$. Let $F=F'-e$.
Clearly, $G\cup F$ contains a subgraph isomorphic to $H$. Let us
assume that $(G-e)\cup F$ contains a subgraph, say $H'$ isomorphic
to $H$. This subgraph is not contained entirely in $G-e$ (as then
$G$ would contain a subgraph isomorphic to $H$ and $G\notin\Phi_H$).
Also, $H'$ is not a subgraph of $F$ ($F$ has one fewer edge than
$H'$). Thus, $\{x,y\}$ is a cutset of $H'$ and $xy$ is not an edge
of $H'$. That implies that $H$ has a 2-element cutset whose elements
are not joined with an edge, a contradiction.

\smallskip
\noindent Parts (4) and (5) of the assertion follow from (1) - (3).
Indeed, if a graph $H$ is 3-connected, then it is 2-connected, too.
Moreover, it vacuously satisfies the requirement that for every
2-element cutset $\{x,y\}$, $xy$ is an edge of $H$. Thus, (4)
follows. If $H$ is a complete graph, then it is a star (if it
consists of only one edge) or a cycle (if it consists of three
edges), or is 3-connected. Thus, (5) follows. \hfill$\Box$

We note that except stars, the $3$-edge path is the only tree $H$
such that $\Phi_H$ is strong. Indeed, if $H$ is a $k$-edge tree
different from a star and a $3$-edge path then we define $G$ to be a
complete graph on $k$ vertices and $e$ to be any edge of $G$. Then,
$G\notin\Phi_H$ (as it has only $k$ vertices). However, for every
graph $F$, either both $G\cup F$ and $(G-e)\cup F$ contain $H$ or
neither does. Thus, $G\equiv^{\Phi_H}_s G-e$, which implies that
$\Phi_H$ is not strong.

If $H$ is a $3$-edge path then $G\not\in\Phi_H$ if and only if every
component of $G$ is a star or a triangle. Let $e=xy$ be an edge in
$G$. If $e$ itself is a component of $G$ then let $F$ be a $2$-edge
path $yzu$, where $z$ and $u$ are new vertices (not in $V(G)$). If
$e$ is an edge of a star centered at $x$ then we define $F$ to be an
edge $yz$, where $z$ is a new vertex. Finally, if $e$ is an edge of
a triangle then let $F$ be an edge $zu$, where $z$ is the third
vertex of the triangle and $u$ is a new vertex. In each case
$G\cup F$ contains a $3$-edge path while $(G-e)\cup F$ does not.
Hence $G\not\equiv_s^{\Phi_H}G-e$ which, by Theorem
\ref{thm:exchange}, implies that $\Phi_H$ is strong when $H$ is the
$3$-edge path.


Theorem \ref{thm:pm} states that for every 3-connected graph $H$,
the property $\Phi_H$ is strong. However, the problem of
characterizing graphs $H$ of connectivity $1$ and $2$ such that
$\Phi_H$ is not strong is open.

\section{Colorability, Edge Colorability and Connectivity}

In the rest of the paper, we discuss strengthening of equivalence
relations arising in the context of some well studied graph-theoretic
concepts: connectivity, colorability and edge colorability.
We start with
a simple lemma.

\begin{lemma}\label{skladowe}
Let $G$ and $H$ be graphs. If $V(G)=V(H)$ and the families of vertex
sets of components of $G$ and $H$ are not the same then there exists
a pair of vertices which are joined by an edge in one of the graphs
$G$ or $H$ and are in two different components in the other graph.
\end{lemma}
{\bf Proof.} By our assumptions there exists a component, say $G'$
in $G$ such that $V(G')$ is not the vertex set of any component in
$H$. Let $H'$ be a component in $H$ such that $V(H')\cap
V(G')\not=\emptyset$ and let $x$ be any element of $V(H')\cap
V(G')$. As $V(G')\not=V(H')$, $V(G')-V(H')\not=\emptyset$ or
$V(H')-V(G')\not=\emptyset$. We will consider the former case only
because both cases are very similar. Let $y\in V(G')-V(H')$. Since
both $x$ and $y$ belong to the component $G'$, there exists a path
$P$ in $G$ joining $x$ with $y$. Clearly, $x$ and $y$ belong to two
different components in $H$ so there is an edge in the path $P$ that
joins vertices that belong to two different components of $H$.
$\hfill\Box$

\subsection{Colorability}\label{sub41}

Let $k$ be a positive integer and let the function $cl$ assign to
every graph the set of its good $k$-colorings (to simplify the
notation, we drop the reference to $k$). We will show that
$\equiv^{cl}_s\ =\ \cong^{cl}\ =\ \cong^{cl}_s$.

\begin{theorem}\label{color1}
Let $k$ be a positive integer. For every graphs $G$ and $H$, the following
conditions are equivalent:
\begin{description}
\item[{\rm (i)}] $G\equiv^{cl}_s H$,
\item[{\rm (ii)}] $G\cong^{cl}H$,
\item[{\rm (iii)}] $G\cong^{cl}_s H$.
\end{description}
\end{theorem}
{\bf Proof.} (i) $\Rightarrow$ (ii). Let $G\equiv^{cl}_s H$. By the
inclusion $\equiv^{cl}_s \ \subseteq\ \equiv^{cl}$, either both $G$ and $H$
are well $k$-colorable or both $G$ and $H$ are not well
$k$-colorable. In the latter case the sets of good $k$-colorings of
$G$ and $H$ are empty so they are equal. Thus assume that both $G$
and $H$ are well $k$-colorable.

Let us suppose $V(G)\not=V(H)$ and assume without loss of generality that
there exists a vertex $v\in V(G)-V(H)$. Denote by $x$ a neighbor of
$v$ in $G$. We define a graph $F$ whose vertex set consists of the
vertices $v, x$ and some $k-1$ vertices which are not in $V(G)\cup
V(H)$. The edges of $F$ join each pair of vertices except $v$ and
$x$. The graph $G\cup F$ contains a complete subgraph $K_{k+1}$ on
the vertex set $V(F)$, so $G\cup F$ is not well $k$-colorable. On
the other hand the graph $H\cup F$ is well $k$-colorable because
both $H$ and $F$ are well $k$-colorable and their only possible
common vertex is $x$. Hence $G\not\equiv^{cl}_s H$, a contradiction.

We have shown that $V(G)=V(H)$. Let us suppose that the sets of good
$k$-colorings of $G$ and $H$ are not equal and assume without loss
of generality that there is a good $k$-coloring ${\cal C}=\{
C_1,C_2,\ldots,C_k\}$ of $G$, which is not a good $k$-coloring of
$H$. We define $F$ to be the complete $k$-partite graph on the set
of vertices $V(G)$, whose monochromatic sets of vertices are
$C_1,C_2,\ldots,C_k$. Clearly, $\cal C$ is the only good
$k$-coloring of $F$. Since $G\cup F=F$, $G\cup F$ has a good
$k$-coloring. On the other hand, $\cal C$ is not a good $k$-coloring
of $H$. Thus, it is not a good $k$-coloring of $H\cup F$, either. We
have shown that $G\not\equiv^{cl}_s H$. This contradiction proves
that the sets of good $k$-colorings of $G$ and $H$ are the same.

\smallskip
\noindent (ii) $\Rightarrow$ (iii) Let us assume that the sets of
good $k$-colorings of $G$ and $H$ are equal. If both these sets are
empty then they remain empty (and so, equal) for the graphs $G\cup
F$ and $H\cup F$, for every graph $F$. If the sets of good
$k$-colorings of $G$ and $H$ are equal but not empty then
$V(G)=V(H)$, as good $k$-colorings of a graph are partitions of the
set of vertices of this graph. Let $F$ be any graph. If both $G\cup
F$ and $H\cup F$ are not well $k$-colorable then $G\cup F\cong^{cl}
H\cup F$. Let now ${\cal C}=\{ C_1,C_2,\ldots,C_k\}$ be a good
$k$-coloring of $G\cup F$. Obviously, ${\cal C}'=\{ C_1\cap
V(G),C_2\cap V(G),\ldots,C_k\cap V(G)\}$ is a good $k$-coloring of
$G$, so by our assumption, ${\cal C}'$ is a good $k$-coloring of
$H$. Thus no block of $\cal C$ contains an edge of $H$ (as
$V(G)=V(H)$). Similarly, no block of $\cal C$ contains an edge of
$F$ because ${\cal C}''=\{ C_1\cap V(F),C_2\cap V(F),\ldots,C_k\cap
V(F)\}$ is a good $k$-coloring of $F$. It follows that $\cal C$ is a
good $k$-coloring of $H\cup F$. Thus the set of good $k$-colorings
of $G\cup F$ is a subset of the set of good $k$-colorings of $H\cup
F$. The converse inclusion can be shown exactly the same way.
Consequently, (iii) holds.

\smallskip
\noindent (iii) $\Rightarrow$ (i) This implication follows from the
obvious inclusion $\cong^{cl}_s \ \subseteq\ \equiv^{cl}_s $ (we
refer to Proposition \ref{prop:2}(1) and note that $\cong^{cl}\
\subseteq\ \equiv^{cl}$). $\hfill\Box$

\begin{theorem}
Let $cl$ be the function assigning to every graph the set of its
good $k$-colorings. If $k\geq 3$ then the problem of deciding if
$G\not\equiv^{cl}_s H$ is NP-complete.
\end{theorem}
{\bf Proof.} By Theorem \ref{color1}, to demonstrate that
$G\not\equiv^{cl}_s H$, it suffices to show a $k$-coloring of one of
the graphs $G$ or $H$ which is not a $k$-coloring of the other
graph. Hence the problem is in NP.

We will now reduce the NP-complete problem of existence of a
$k$-coloring of a graph (cf. Garey and Johnson \cite{GJ}) to our
problem. Let $G'$ be a graph. We define $G$ to be the graph obtained
from $G'$ by adding an edge $xy$, where $x$ and $y$ are two new vertices.
We define $H$ to be the graph obtained from $G$ by adding an edge $zx$,
where $z$ is some vertex of $G'$. We will prove that $G'$ is $k$-colorable
if and only if $G\not\equiv^{cl}_s H$.

Let us suppose first that $G'$ is $k$-colorable. Then there exists a
$k$-coloring of $G$ such that both vertices $x$ and $z$ belong to
the same block of the $k$-coloring. This $k$-coloring of $G$ is not
a $k$-coloring of $H$ because $xz$ is an edge in $H$. Thus
$G\not\cong^{cl} H$ and, by Theorem \ref{color1}, $G\not\equiv^{cl}_s H$.
Conversely, if $G'$ is not $k$-colorable then neither is $G$ nor
$H$. Hence, the sets of $k$-colorings of both $G$ and $H$ are empty and,
consequently, equal. By Theorem \ref{color1} again $G\equiv^{cl}_s H$.
$\hfill\Box$

\begin{remark}
For $k=1$, that is, when the function $cl$ assigns to a graph the
set of its good $1$-colorings, all four equivalence relations
$\equiv^{cl},\cong^{cl},\equiv^{cl}_s $ and $\cong^{cl}_s $ coincide
and for every pair of graphs $G$ and $H$, $G\equiv^{cl} H$ if and
only if $G=H=\emptyset$ or $G\not=\emptyset\not=H$.
\hfill$\Box$
\end{remark}

For $k=2$ the problem if $G\not\equiv^{cl}_s H$ is solvable in polynomial
time. It is a consequence of the following fact.

\begin{proposition}
Let $cl$ be the function assigning to every graph the set of its good
$2$-colorings. For any two graphs $G$ and $H$, $G\equiv^{cl}_s H$ if
and only if either none of the graphs $G$ and $H$ is bipartite or
\begin{enumerate}
\item both $G$ and $H$ are bipartite
\item $V(G)=V(H)$
\item $G$ and $H$ have the same families of vertex sets of connected
components, and
\item connected components with the same vertex sets in $G$ and $H$ have
the same bipartitions.
\end{enumerate}
\end{proposition}
{\bf Proof.} $(\Ra)$ Let $G\equiv^{cl}_s H$ and assume that at least
one of the
graphs $G$ or $H$ is bipartite. Otherwise the necessity holds. As
$\equiv^{cl}_s \ \subseteq\ \equiv^{cl}$, both $G$ and $H$ are bipartite.
Hence the sets of good $2$-colorings of $G$ and $H$ are nonempty and
they are equal by Theorem \ref{color1}. Since good $k$-colorings in
a graph are partitions of the vertex set of this graph,
$V(G)=V(H)$.

Let us suppose the families of vertex sets of components of $G$ and
$H$ are not the same. Then, by Lemma \ref{skladowe}, there exists a
pair of vertices that are joined by an edge in one of the graphs $G$
or $H$ and are in different components in the other graph. Without
loss of generality we can assume that there are vertices, say $x$
and $y$, that are joined by an edge in $G$ but belong to two
different components in $H$. Let us denote by $H_x$ (respectively,
$H_y$) the two components in $H$ that contain $x$ (respectively, $y$)
and by $V_x$ (respectively, $V_y$) the monochromatic class of $H_x$
(respectively, $H_y$) that contains $x$ (respectively, $y$). Let $F$
be the complete bipartite graph on the set of vertices $V(H_x)\cup
V(H_y)$ whose monochromatic classes are $V_x\cup V_y$ and the other
one $(V(H_x)\cup V(H_y))-(V_x\cup V_y)$. The graph $H\cup F$ is
bipartite while $G\cup F$ is not because it contains the graph $F+
xy$ which has an odd cycle. By the definition of the relation
$\equiv^{cl}_s $, $G\not\equiv^{cl}_s H$, a contradiction. Hence for
each component in $G$, there is a component in $H$ with the same
vertex set. By symmetry, $G$ and $H$ have the same families of
vertex sets of their connected components.

Let us suppose now that for some two components $G'$ of $G$ and $H'$
of $H$ with the same vertex sets, the bipartitions of $G'$ and $H'$
are not the same. Then there exists a pair of vertices $u$ and $v$
such that they both are in the same monochromatic class in $G'$ but
in the different monochromatic class in $H'$. The graph $(G\cup H')+
uv$ is not bipartite because there is a path of an even length
joining $u$ and $v$ in $G$ so $(G\cup H')+ uv$ contains an odd
cycle. On the other hand the graph $(H\cup H')+ uv=H+uv$ is
bipartite. Hence $G\not\equiv^{cl}_s H$, a contradiction. This
completes the proof of necessity.

\smallskip
\noindent
$(\La)$ If both $G$ and $H$ are not bipartite then for every graph $F$, both
$G\cup F$ and $H\cup F$ are not bipartite so $G\equiv^{cl}_s H$. Let
now both $G$ and $H$ be bipartite.
Let us denote by $\cal C$ a good $2$-coloring of $G$ and suppose $\cal C$
is not a good $2$-coloring of $H$. Then there exists an edge, say
$e$, in $H$ whose both ends are contained in the same block of $\cal
C$. Let $H'$ be the component of $H$ that contains this edge. By our
assumptions, there is a component $G'$ of $G$ that has the same
vertex set and the same bipartition as $H'$. Thus the edge $e$ has
its ends in two different blocks of this bipartition of $G'$ so in
two different blocks of $\cal C$ as well. This contradiction proves
that every good $2$-coloring of $G$ is a good $2$-coloring of $H$.
In a very similar way one can prove that every good $2$-coloring of
$H$ is a good $2$-coloring of $G$. By Theorem \ref{color1}, we
conclude that $G\equiv^{cl}_s H$. $\hfill\Box$

\subsection{Edge 2-colorability}

We will now consider the property of edge 2-colorability. We note that
a graph is \emph{edge 2-colorable} if and only if each of its connected
components is a path or a cycle of an even length We denote by
$\equiv^{e2c}$ the equivalence relation in which two graphs are
equivalent if and only if both are edge 2-colorable or neither of
the two is.

The following theorem characterizes the relation $\equiv^{e2c}_s$.

\begin{theorem}
Let $G$ and $H$ be graphs. Then $G\equiv^{e2c}_s H$ if and only if
at least one of the following conditions holds
\begin{enumerate}
\item Both $G$ and $H$ are not edge 2-colorable (each contains an odd cycle
or a vertex of degree at least 3)
\item Both $G$ and $H$ are edge 2-colorable (no odd cycles and maximum
degree at most 2), $V(G)=V(H)$, and for every even (odd) path component
in $G$ there is an even (odd) path component in $H$ with the same endpoints.
\end{enumerate}
\end{theorem}
{\bf Proof.} $(\Ra)$ If both $G$ and $H$ are not edge 2-colorable,
there is nothing left to prove. Since $G\equiv^{e2c}_s H$, then
$G\equiv^{e2c}H$ and, consequently, both $G$ and $H$ are edge
2-colorable. Let us suppose that $V(G)\not=V(H)$, say
$V(H)-V(G)\not=\emptyset$. Let $u\in V(H)-V(G)$. We denote by $v$
and $w$ some new vertices (occurring neither in $G$ nor in $H$). We
define $F=\{vu,wu\}$. Clearly, $G\cup F$ is edge 2-colorable, while
$H\cup F$ is not. Thus, $G\not \equiv^{e2c}_s H$, a contradiction,
so $V(G)=V(H)$.

Suppose now there is a vertex $u$ of $G$ such that $deg_G(u)=1$ and
$deg_H(u)=2$. Let $v$ be a new vertex (occurring neither in $G$ nor
in $H$). We define $F=\{vu\}$. Obviously, as before, $G\cup F$ is
edge 2-colorable, while $H\cup F$ is not. Hence, $G\not
\equiv^{e2c}_s H$, a contradiction. By the symmetry argument it
follows that for every vertex $u$, $deg_G(u)=deg_H(u)$.

Let now $P$ be a path in $G$ with endpoints $a$ and $b$. By the
property proved above, $deg_H(a)=deg_H(b)=1$. Let us assume that $a$
and $b$ are endpoints of two different paths in $H$. We select a new
vertex, say $u$. If $P$ has odd length, we define $F=\{au,ub\}$.
Otherwise, we define $F=\{ab\}$. Clearly, $G\cup F$ contains an odd
cycle and so, it is not edge 2-colorable. On the other hand, $H\cup
F$ does not contain any odd cycles and so, it is edge 2-colorable, a
contradiction.

Thus, $a$ and $b$ are the endpoints of the same path in $H$, say
$P'$. It remains to prove that the length of $P'$ is of the same
parity as the length of $P$. If $G+ab$ is edge 2-colorable, then the
cycle $P+ab$ has even length. As $G\equiv^{e2c}_s H$, $H+ab$ is edge
$2$-colorable too. Thus, the cycle $P'+ab$ is also even, so both
paths $P$ and $P'$ are of odd length. Similarly, if $G+ab$ is not
edge 2-colorable, then $P+ab$ is an odd cycle. It follows that
$P'+ab$ is an odd cycle, too and, consequently, $P$ and $P'$ are
both of even length.

\smallskip
\noindent $(\La)$ Let $F$ be graph and let us assume that $G\cup F$
is edge 2-colorable. It follows that $G$ is edge 2-colorable and so,
$H$ is edge 2-colorable, too. Moreover, no vertex in $G\cup F$ has
degree 3 or more. By our assumptions, the same holds for $H\cup F$
because the degrees of vertices in $G$ and $H$ are the same. Let us
consider any cycle $C$ in $H\cup F$. If $C\cap F=\emptyset$, then
$C\subseteq H$. Consequently, $C$ is even. Thus, let $F'=C\cap F$.
By our assumptions, adding $F'$ to $G$ results in exactly one new
cycle in $G$, say $C'$. Moreover, the parity of the lengths of $C'$
and $C$ is the same. Since $G\cup F$ is edge 2-colorable and
contains $C'$, $C'$ is even. Thus, $C$ is even, too. It follows that
$H\cup F$ is edge 2-colorable. By symmetry, for every graph $F$,
$G\cup F$ is edge 2-colorable if and only if $H\cup F$ is edge
2-colorable. \hfill$\Box$

\subsection{Connectivity}

By a {\em cutset} in a connected graph we mean a set of vertices in
this graph whose deletion makes this graph disconnected. A set $C$
of vertices in a disconnected graph $G$ is a {\em cutset} of $G$ if
$C= \emptyset$ or, for some component $G'$ of $G$, $C\cap V(G')$ is
a cutset of $G'$. Clearly, $C\not=\emptyset$ is a cutset of $G$ if
and only if it separates some pair of vertices in $G$. Let us
observe that the only graphs without any cutsets are the complete
graphs.

Let $cs$ be a function that assigns to every graph the set of its
cutsets of cardinality smaller than $k$, where $k\geq 1$ (as in the case of
colorability, to simplify the notation, we drop the reference to
$k$). Thus, $G\equiv^{cs} H$, if either both graphs $G$ and $H$ have
a cutset of cardinality less than $k$ or both do not have such a
cutset. We shall characterize now the relation $\equiv^{cs}_s$.

\begin{lemma}\label{spoj1}
If $G\equiv^{cs}_s H$ then $V(G)=V(H)$.
\end{lemma}
{\bf Proof.} Let us suppose that $V(G)\not=V(H)$. We can assume without
loss of generality that there exists a vertex $x\in V(H)-V(G)$.

We will first assume that $k=1$. Let $z$ be any vertex in $H$ different
from $x$ and let $y$ be a new vertex (not in $H$ nor $G$). We define
$F=\{zu\st u\in V(H)-\{x,z\}\}\cup\{xy\}$. It follows that $H\cup F$ is
connected (that is, has no cutsets of size 0) and $G\cup F$ is not
connected (the edge $xy$ is separated from the rest of the graph).
Thus, $G\not\equiv^{cs}_s H$, a contradiction.


Thus, from now on, we assume that $k\geq 2$. Let
$K$ 
be a set of $k$ vertices which are not in $V(G)\cup V(H)$
and let $\ell=\max(k-\mbox{deg}_H(x),1)$. Since $x$ is not an isolated
vertex in $H$ and $k\geq 2$, $\ell<k$. We define $F$ to be the graph
obtained from
the complete graph on $(V(H)\cup K)-\{ x\}$ by adding the vertex $x$
and edges joining $x$ with all vertices of some $\ell$-element
subset $L$ of $K$. The graph $H\cup F$ can be obtained from $F$ by
adding the edges incident in $H$ with $x$. As $\mbox{deg}_{H\cup
F}(x)=\mbox{deg}_H(x) + \ell\geq k-\ell+\ell=k$, the graph $H\cup F$
does not have a cutset of cardinality smaller than $k$.
On the other hand,
the set $L$ 
is a cutset in $G\cup F$. Indeed, $L$ is a cutset of the
component of $G\cup F$ containing $K\cup\{ x\}$ as it separates $x$
from the rest of the graph. Since $|L|<k$, $G\cup F$ has a cutset of
cardinality smaller than $k$. Thus, $G\not\equiv^{cs}_s H$, a
contradiction, and so $V(G)=V(H)$ follows. $\hfill\Box$

\begin{theorem}\label{spoj2}
Let $G$ and $H$ be graphs. Then, $G\equiv^{cs}_s H$ if and only if
$V(G)=V(H)$ and for every set $C\subseteq V(G)$ such that $|C|<k$,
the families of vertex sets of components of $G-C$ and $H-C$ are the
same.
\end{theorem}
{\bf Proof.} $(\La)$ To show sufficiency assume that $G\not\equiv^{cs}_s H$.
Then, for some graph $F$, $G\cup F\not\equiv^{cs} H\cup F$. We can
assume without loss of generality that $G\cup F$ has a cutset $C$ of
cardinality smaller than $k$ and $H\cup F$ does not have such a
cutset. As $(G\cup F)-C$ is disconnected, there are vertices $x$ and
$y$ which belong to two different components of $(G\cup F)-C$. On
the other hand $(H\cup F)-C$ is connected so there exists a path,
say $P$, joining $x$ and $y$ in $(H\cup F)-C$. Let $e$ be any edge
in $P$ which does not belong to $F$. Then, clearly, $e$ is an edge
of $H-C$. We denote by $H'$ the component of $H-C$ which contains
the edge $e$. Since the families of vertex sets of components of
$G-C$ and $H-C$ are the same, there is a component of $G$ that
contains the edge $e$. Consequently there exists a path, say $P_e$
in $G-C$ joining the ends of the edge $e$. Let us replace in the
path $P$ every edge $e$ which is not in $F$ by the path $P_e$. The
resulting graph is a connected subgraph of $(G\cup F)-C$ containing
the vertices $x$ and $y$. We have got a contradiction with the
definition of $x$ and $y$. Thus our initial assumption that
$G\not\equiv^{cs}_s H$ was false, so $G\equiv^{cs}_s H$.

\smallskip
\noindent
$(\Ra)$
We pass on to the proof of necessity. By Lemma \ref{spoj1},
$V(G)=V(H)$.

Let $C\subseteq V(G)$, $|C|<k$. Let us suppose the families of
vertex sets of components of $G-C$ and $H-C$ are not the same. By
Lemma \ref{skladowe}, there exists a pair of vertices which are
joined by an edge in one of the graphs $G-C$ or $H-C$ and are in two
different components in the other graph. Without loss of generality
we assume that there is an edge $e$ in $H-C$ whose ends belong to
two different components in $G-C$. Let $G'$ be one of these
components. We denote by $L$ a set of $\ell=k-1-|C|$ vertices which
are not in $V(G)\cup V(H)$. We define $F$ to be the graph on the set
of vertices $V(G)\cup L$ in which every pair of vertices in $V(G')$
is joined by an edge, every pair of vertices in $V(G)-V(G')$ is
joined by an edge and every vertex in $C\cup L$ is joined by an edge
with every other vertex of $F$. Clearly, $C'=C\cup L$ is a cutset in
$G\cup F=F$ of cardinality $|C'|=|C|+k-1-|C|=k-1$. Let us consider
any cutset $C''$ in $H\cup F$. Since $V(H\cup F)=V(H)\cup
L=V(F)=V(G)\cup L=V(G\cup F)$ and every vertex in $C'$ is joined by
an edge with every other vertex of $H\cup F$, $C''\supseteq C'$. Let
us observe that $(H\cup F)-C'\supseteq (F-C')+ e$. The last graph is
a connected spanning subgraph of $(H\cup F)-C'$, so $C'$ is not a
cutset of $H\cup F$. We have shown that every cutset in $H\cup F$
has at least $k$ vertices which shows that $G\not\equiv^{cs}_s H$.
This contradiction completes the proof. $\hfill\Box$

\begin{remark}
\label{rem:3} It follows from Theorem \ref{spoj2} that for every
fixed $k$, the problem to decide if $G\equiv^{cs}_s H$ is polynomial
time solvable, where $cs$ is the function assigning to every graph
the set of its cutsets of cardinality smaller than $k$. It is an
open question what the complexity status of this problem is when $k$
is a part of the instance. \hfill$\Box$
\end{remark}

Next, we will show that $\equiv^{cs}_s\ =\ \cong^{cs}_s$. That is,
$\cong^{cs}$ despite being more precise than $\equiv^{cs}$ has the same
strengthening.

\begin{theorem}\label{conn2}
$G\equiv^{cs}_s H$ if and only if $G\cong^{cs}_s H$.
\end{theorem}
{\bf Proof.} ($\Leftarrow$) This implication follows from a
generally true inclusion $\cong^{f}_s\ \subseteq\ \equiv^{f}_s$.

\smallskip
\noindent ($\Rightarrow$) Let $G$ and $H$ be graphs such that
$G\equiv^{cs}_s H$. We assume that there is a graph $F$ such that
$G\cup F\not\cong^{cs} H\cup F$. Since $G\equiv^{cs}_s H$, by
Proposition \ref{prop:1}(3), $G\cup F\equiv^{cs}_s H\cup F$. By
Lemma \ref{spoj1}, $V(G\cup F)=V(H\cup F)$. As $G\cup
F\not\cong^{cs} H\cup F$, we can assume without loss of generality
that there exists a cutset $C$ in $G\cup F$, $|C|<k$, which is not a
cutset in $H\cup F$. Let $x$ and $y$ be a pair of vertices which
belong to two different components of $(G\cup F)-C$. We denote by
$G'$ the component of $(G\cup F)-C$ that contains $x$. Clearly,
there is a path in $(H\cup F)-C$ that joins the vertices $x$ and
$y$. Consequently there exists an edge, say $e$, in $(H\cup F)-C$
with one vertex in $V(G')$ and the other one in $V((H\cup
F)-C)-V(G')=V((G\cup F)-C)-V(G')$.

Let $L$ be a set of cardinality $k-1-|C|$ of vertices not in
$V(G\cup F)$. We define $K$ to be the graph on the set of vertices
$V(G\cup F)\cup L$ in which every pair of vertices in $V(G')$ is
joined by an edge, every pair of vertices in $V((G\cup F)-C)-V(G')$
is joined by an edge and every vertex in $C'=C\cup L$ is joined by
an edge with every other vertex of $K$. Clearly, $C'$ is a cutset in
$G\cup F\cup K=K$ of cardinality $|C'|=|C|+k-1-|C|=k-1$. Since
$G\equiv^{cs}_s H$, the graph $H\cup F\cup K$ has a cutset $C''$ of
cardinality smaller than $k$. We observe that $C''\supseteq C'$
because $V(H\cup F\cup K)=V(K)$ and every vertex in $C'$ is joined
by an edge with every other vertex in $V(H\cup F\cup K)$. As
$|C''|<k$, $C''=C'$. By the definition of $K$, $V((H\cup F\cup
K)-C')=V(K-C')=V((G\cup F)-C)$ and the graph $(H\cup F\cup K)-C'$
contains complete graphs on the sets of vertices $V(G')$ and
$V((G\cup F)-C)-V(G')$ as subgraphs. The graph $(H\cup F\cup K)-C'$
contains the edge $e$ whose one end is in one of the complete
subgraphs mentioned above and the other one in the other complete
subgraph. Therefore the graph $(H\cup F\cup K)-C'$ is connected so
$C''=C'$ is not a cutset in $H\cup F\cup K$, a contradiction.

We have shown that $G\cup F\cong^{cs} H\cup F$ for every graph $F$, so
$G\cong^{cs}_s H$. $\hfill\Box$

\begin{remark}
Unlike in the case of $k$-coloring, $\cong^{cs}\ \not=\
\cong^{cs}_s$. A simple example for $k=2$ is shown in Figure
\ref{conn}.

\begin{figure}
\centerline{\includegraphics[scale=0.40]{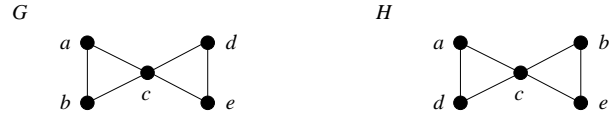}}
\caption{The set $\{c\}$ is the only one-element cutset of $G$ and of
$H$. Thus,
$G\cong^{cs} H$. However, components of $G-c$ and $H-c$ are different and so,
$G\not\cong^{cs}_s H$.\hfill$\Box$}
\label{conn}
\end{figure}
\end{remark}

\vspace*{-0.2in} For a given graph $G$ it would be interesting to
find a smallest (with respect to the number of edges) subgraph $G'$
of $G$ such that $G'\equiv^{cs}_s G$. We observe, however, that even
for $k=2$, it is a difficult problem. Indeed, let us consider the
problem of deciding if for a given graph $G$ and an integer $m$
there exists a subgraph $G'$ of $G$ such that $G\equiv^{cs}_s G'$
and $G'$ has at most $m$ edges. One can easily verify that the
problem is in the class NP (cf. Remark \ref{rem:3}). Moreover it is
NP-complete because it follows from Theorem \ref{spoj2} that for a
$2$-connected graph $G$ and $m=|V(G)|$ the problem asks for the
existence of a hamiltonian cycle in $G$.

Let $\Phi$ be the set of all graphs with a cutset of cardinality smaller
than $k$. It is obvious that $\equiv^\Phi\ =\ \equiv^{cs}$. Let us denote
by $\Psi$ the set of graphs that are not $k$-connected. One can easily
observe that a graph $G\in\Psi$ if and only if $G$ has a cutset of
cardinality smaller than $k$ or $|V(G)|\leq k$. As the two relations
are closely related, it is natural to ask if $\equiv^\Phi_s\ =\ \equiv^\Psi_s$.
We will answer this question positively.

\begin{theorem}
\label{prop:last}
$G\equiv^\Phi_s H$ if and only if $G\equiv^\Psi_s H$.
\end{theorem}
{\bf Proof.} ($\Rightarrow$) Let us suppose $G\equiv^\Phi_s H$ but
$G\not\equiv^\Psi_s H$. Then, without loss of generality, there exists
a graph $F$ such that $G\cup F\in \Psi$ and $H\cup F\not\in
\Psi$. In this case, $H\cup F$ has no cutsets of cardinality
smaller than $k$ and $|V(H\cup F)|>k$. Since $G\equiv^\Phi_s H$, $G\cup
F$ has no cutsets of cardinality smaller than $k$. Thus, $|V(G\cup
F)|\leq k$ because $G\cup F\in \Psi$. Hence $|V(G\cup
F)|<|V(H\cup F)|$, a contradiction because, by Lemma \ref{spoj1},
$V(G)=V(H)$, so $V(G\cup F)=V(H\cup F)$, as well.

\smallskip
\noindent
($\Leftarrow$) Let us suppose now that $G\equiv^\Psi_s H$ but
$G\not\equiv^\Phi_s
H$. Then, without loss of generality, there exists a graph $F$ such
that $G\cup F$ has a cutset of cardinality smaller than $k$ and
$H\cup F$ does not have such a cutset. It follows that
$G\cup F\in \Psi$. Since $G\equiv^\Psi_s H$, $H\cup F\in \Psi$.
Thus, $|V(H\cup F)|\leq k$ and, consequently, $H\cup F$ is a complete
graph on at most $k$ vertices. As $G\equiv^\Psi_s H$, $G\cup H\cup
F\equiv^\Psi H\cup H\cup F=H\cup F$, so $G\cup H\cup F\in \Psi$.

Let us suppose first that $G\cup H\cup F$ has a cutest $C$ of cardinality
smaller that $k$. Let $G'$ be one of the components of $(G\cup H\cup
F)-C$ and let $G''=(G\cup H\cup F)-C-G'$. The complete graph $H\cup
F$ has common vertices with at most one of the graphs $G'$ and
$G''$, say with $G''$. Let $H'$ be a complete graph on $k+1$
vertices that contains $H\cup F$ and has no common vertices with
$G'$. Now, $G\cup H\cup F\cup H'\in \Psi$ because $C$ is its
cutset, while $H\cup H\cup F\cup H'=H'\not\in\Psi$ because it is a
complete graph on $k+1$ vertices. Hence $G\not\equiv^\Psi_s H$, a
contradiction.

We have proved that $G\cup H\cup F$ has no cutsets of cardinality
smaller that $k$ so, as $G\cup H\cup F\in \Psi$, $|V(G\cup H\cup
F)|\leq k$. Consequently, $|V(G\cup H)|\leq k$. Since
$G\not\equiv^\Phi_s H$, $G\not=H$ have different edge sets. We can
assume without loss of generality that there is an edge, say $e=xy$
in $H$, which is not an edge in $G$. Let $K$ be a graph obtained
from the complete graph on a $(k+1)$-element set of vertices
containing $V(G\cup H)$ by deleting the edge $e$. Clearly, $G\cup
K=K$ but $H\cup K$ is the complete graph on $k+1$ vertices. The
former graph has a cutset $V(K)-\{ x,y\}$ of cardinality $k-1$ and
the latter graph has no cutsets. Hence $G\cup K\in\Psi$ and $H\cup
K\not\in\Psi$ so $G\not\equiv^\Psi_s H$. This contradiction proves
that $G\equiv^\Phi_s H$. $\hfill\Box$

\section{Open problems and further research directions}

We do not know of any past research concerning the strengthening of
equivalence relations on graphs.
Nevertheless, it seems to us that this is a natural concept worth
further investigations. In this paper, we focused on the strengthening
of the equivalence relations that are determined by graph properties.
For many properties $\Phi$ we studied, the relations $\equiv^\Phi_s$
turned out to have a very simple structure (for instance, they broke both
equivalence classes of $\equiv^\Phi$, $\Phi$ and $\overline{\Phi}$, into
singletons, that is, were identities; or kept $\Phi$ as an equivalence
class and broke the other one into singletons). In several cases, however,
(for the properties $\Phi$ of being $k$-connected, $k$-colorable, and edge
$2$-colorable), the structure of the relations $\equiv_s^\Phi$ turned out
to be more complex and so more interesting, too. Therefore, a
promising research direction could be to identify and study
additional natural graph properties $\Phi$, for which the relations
$\equiv_s^\Phi$ have a nontrivial structure. Establishing characterizations
of the relations $\equiv_s^\Phi$, and determining the complexity of
deciding whether for two given graphs $G$ and $H$, $G\equiv^\Phi_s H$ holds,
are particularly interesting and important. Given the results of our paper,
it seems that the properties of being edge $k$-colorable and edge
$k$-connected are natural candidates for this kind of investigations.
In the former case we were only able to find a characterization of the
relation $\equiv_s^\Phi$, when $\Phi$ is the property of being edge
$2$-colorable. The theorem we proved in this case suggests that for
an arbitrary $k$ the structure of the relation $\equiv_s^\Phi$ may be quite
complicated, which makes the problem a challenge. In the latter
case, we feel there may be strong similarities with the strengthening of the
property of $k$-connectivity but do not have any specific results.

There are a few open problems directly related to the results of
this paper. One of them is to establish the computational complexity
of the problem to decide if $G\equiv^\Phi_s H$ (given graphs $G$ and $H$,
and an integer $k$), when $\Phi$ is the property of being
$k$-connected. For a fixed $k$, it follows from our Theorem
\ref{spoj2} that the problem is solvable in polynomial time. The
question is open however, when $k$ is a part of the instance.

Another problem concerns the property $\Phi_H$ of containing a
subgraph isomorphic to $H$. It was shown in Theorem \ref{thm:pm}
that the relation $\equiv^{\Phi_H}_s$ has some very simple structure
for many graphs $H$. The question arises what is the structure of
$\equiv^{\Phi_H}_s$ for all other graphs $H$.

There are also several natural general questions concerning the
concept of strengthening of an equivalence relation in graphs. For
example it would be interesting to find a general condition for the
function $f$ that guarantees that the relations $\equiv^f_s$ and
$\cong^f_s$ are equal.

Another general problem is to establish conditions that ensure that
there exists a weakest equivalence relation $\equiv'$ such that
$\equiv'_s$ is the same as $\equiv_s$ and, whenever it is so, to
find this $\equiv'$. In some cases, the problem is easy. For
example, the relations $\equiv^{cl}_s$ and $\cong^{cl}_s$ studied in
Subsection \ref{sub41} are equal (see Theorem \ref{color1}) and the
relation $\equiv^{cl}$ is strictly weaker than $\cong^{cl}$. As
$\equiv^{cl}$ has only two equivalence classes and the strengthening
of the total relation (the only possible weakening of $\equiv^{cl}$)
is also the total equivalence relation, $\equiv^{cl}$ is the weakest
relation $\equiv'$ such that $\equiv'_s$ and $\cong^{cl}_s$ are
equal. This observation does not generalize to other properties. For
instance, the relation $\equiv^{cs}$ is not the weakest relation
$\equiv'$ such that $\equiv'_s$ and $\cong^{cs}_s$ are the same.
Indeed, the strengthening of $\equiv^\Psi$ (where $\Psi$ is the
property considered in Theorem \ref{prop:last}) coincides with
$\cong^{cs}_s$, yet $\equiv^{cs}$ is incomparable with $\equiv^\Psi$
(and so, in this case, there is no weakest relation of the desired
property). Thus, in general, the two problems mentioned above seem
to be nontrivial.




\end{document}